\documentclass[12pt,english]{article}
\usepackage[cp1251]{inputenc}
\usepackage{geometry}
\usepackage{color}
\usepackage{amsmath}
\usepackage{amssymb}
\usepackage{esint}
\usepackage{amsthm}

\makeatletter
\newtheorem{remark}{Remark}

\makeatother

\usepackage{babel}
\newtheorem{theorem}{Theorem}
\newtheorem{lemma}{Lemma}

\begin{document}
\global\long\def\E{\mathbb{E}}
\global\long\def\P{\mathbb{P}}
\global\long\def\N{\mathbb{N}}
\global\long\def\ind {\mathbb{I}}

\title{{\normalsize\tt\hfill\jobname.tex}\\
On convergence of 1D Markov diffusions to heavy-tailed  invariant density}

\author{O.A. Manita\footnote{Moscow State University, Moscow, Russia; 
email: oxana.manita @ gmail.com}, 
A.Yu. Veretennikov\footnote{University of Leeds, UK, \& National Research University Higher School of Economics, \& Institute for Information Transmission Problems, Moscow, Russia; email: a.veretennikov @ leeds.ac.uk. For this author the work has been funded by the Russian Academic Excellence Project '5-100' (the sections 1 -- 2, the setting in the section 3  and both Lemmata) and by the Russain Science Foundation project no. 17-11-01098 (steps 1, 3, 8 -- 9 of the proof of the Theorem 1).}
}
\maketitle
\begin{abstract}
Rate of convergence  is studied for a diffusion process on the half line with a non-sticky reflection to a heavy-tailed 1D invariant distribution which density on the half line has a polynomial decay at  infinity. Starting from a standard receipt which guarantees some polynomial convergence, it is shown how to construct
a new non-degenerate diffusion process on the half line which converges to the same invariant
measure exponentially fast uniformly with respect to the initial data. 
\end{abstract}


\section{Introduction}

A topical area of Markov Chain Monte Carlo (MCMC) in theoretical statistics is around the 
following problem: given a fixed ``target'' density or distribution known up to a constant multiplier -- 
a normalizing constant -- how to construct a (Markov) process which would have this density as a (unique) 
invariant one and which would converge to this invariant one with a rate that could be theoretically 
evaluated? In particular, a permanent great interest in recent decades was about dealing with 
``heavy-tailed'' densities with a polynomial decay at infinity. 
With this problem in mind, let us consider a polynomially decreasing probability density $\pi$ on 
the line \(\mathbb R^1\); in the precise setting it will be restricted to the half-line \(\mathbb R^1_+\). The question under consideration in this paper is constructing a Markov 
diffusion process with invariant measure $\pi(x)dx$ such that this measure is invariant for the constructed 
process and, moreover, so that an exponential convergence in total variation to the invariant 
distribution holds. 

This problem has certain deep relations to ergodicity and to the Perron -- Frobenius theorem for  Markov 
chains with finite state space, to spectral gap for semigroup generators, to upper and  
lower bounds for convergence to stationarity; yet, a spectral gap in this paper is not used. The literature in this area is huge and we only mention a few important references related to the subject of the paper more or less directly (see \cite{AitSahalia}, 
\cite{Cattiaux}, \cite{Fort},  \cite{Kovchegov}, \cite{kulik-leonenko}, \cite{Eva}, 
\cite{MenshPopov},  \cite{MenshPopov2}, et al.; also, see further references theiren).

The paper consists of four sections, the first one being this Introduction. In the section 2 
two known receipts of constructing an SDE with a given stationary measure are shown: one is an SDE with a unit diffusion coefficient while  another one is an 
SDE with an affine  drift. In the section 3 we state the main result of this paper, and in the section 
4 its proof is provided. The construction is based on the first one of the standard receipts from the 
section 2 and on a random time change. The proof uses certain recurrence type hitting time moment bounds introduced earlier in \cite{ayv_grad_drift}.

\section{Quick review: two standard receipts on $\mathbb R^1$}
\subsection{Receipt 1: SDE with a unit diffusion coefficient}

Suppose a continuous and differentiable strictly positive probability density $\pi$
on $\mathbb{R}^{1}$  decreases at infinity polynomially,
i.e. there exist constants $c>0$ and $m>1$ such that for any $x$, 
\begin{equation}
	c\left(1+\left|x\right|\right)^{-m}\leq
		\pi\left(x\right)\leq 
			c^{-1}\left(1+\left|x\right|\right)^{-m}.
\label{eq:dens}
\end{equation}
Here $m>1$ is required so that the function $\pi$ were integrable; for further claims a bit more restrictive condition  \(m>3\) will be assumed in the sequel. 

On a probability space $\left(\Omega,\mathcal{F},\mathbb{P}\right)$ let us 
fix a standard Wiener process $W_t$ with its natural filtration  $\mathcal{F}_{t}=\mathcal{F}_{t}^W$ (as usual, $\mathbb P$ -- completed). 
On this probability space 
consider a Langevin diffusion $Y_t$  given by an SDE
\begin{equation}\label{eq:lang}
dY_{t}=dW_{t}+b\left(Y_{t}\right)dt, \quad Y_0=\xi,
\end{equation}
with an arbitrary nonrandom initial value \(\xi\), where 
\begin{equation}\label{eqb}
b(x)=\frac{1}{2}(\ln\pi(x))'.
\end{equation}
If there is no explosion then this equation possesses a strong solution \cite{Ve81}.
Random initial values will be mentioned briefly in the section \ref{sec22} and in principle could have been allowed here, too.  It is assumed that two derivatives \(\pi'\) and $\pi''$ exist and that the drift \(b\) is locally bounded; 
its global boundedness is not required because, as it turns out \cite{ayv_grad_drift2, ayv_grad_drift}, a no blow-up is guaranteed just by the  
assumption  (\ref{eq:dens}) on the function \(\pi\) only (see below the details). Emphasize that despite the assumed two derivatives, the only quantitative assumption will be just on $\pi$ itself given in (\ref{eq:dens});  also there is a hypothesis that the assumption about  $\pi''$ could be dropped, and this is the reason why we refer to \cite{Ve81} instead of more standard results under a local Lipschitz condition on the drift, while talking on strong solutions in the sequel.

~

The receipt (\ref{eq:lang})--(\ref{eqb}) is,  actually, a continuous time analogue of (one of) a standard 
MCMC receipt(s) in discrete time after a suitable limit. We do not recall it because the paper does not rely upon this limiting procedure; however, this is likely to  signify a possible link to MCMC algorithms in discrete time. Obviously $\pi(x)dx$ is the (unique) invariant distribution of the process $Y_t$; this 
can be shown explicitly by checking the Kolmogorov equation for the invariant distribution. 
The process $Y_t$ is ergodic and has a polynomial rate of convergence to the stationary
distribution with density $\pi$ 
\cite[Theorem 1]{ayv_grad_drift2} (under a bit more restricted conditions 
see also \cite{ayv_grad_drift}), at least, if \(m\) is not too small. Note that unlike in most 
of other works on convergence or mixing rates, Lyapunov functions are not used here, as they were not used in \cite{ayv_grad_drift2, ayv_grad_drift}. 
For close 
results for some particular  distributions and for homogeneous Markov processes under 
various assumptions (usually more restrictive because of explicit assumptions about 
the derivative \(\pi'\)) see also \cite{Abu-Ver09, Abu-Ver09b, kulik-leonenko}; for discrete 
time examples -- that is, actually, about MCMC 
algorithms -- see, e.g., \cite{MP, Fort, MenshPopov2} and further references therein. Emphasize that the 
assumptions 
in \cite[Theorem 1]{ayv_grad_drift2} as well as in \cite[Theorem 1]{ayv_grad_drift} are 
all on $\pi$ and not on~$\pi'$ except that the latter derivative exists and that \(b\) is 
locally bounded. It remains to be our goal to avoid any assumptions on $\pi'$ beyond 
its existence and local boundedness of \(b\) in the sequel.

Moreover, it is known that  under the {\it additional  assumption} about $\pi'$,
\begin{equation}\label{binfty}
\liminf_{|x| \to \infty} x b(x) = 
\liminf_{|x| \to \infty} \frac{x\pi'(x)}{2\pi(x)} = -r  < - 3/2,  
\end{equation}
the {\em beta-mixing rate} of $Y_{t}$  is {\em no faster than polynomial, $\ge C t^{-k}$ with any $k> r-1/2$ and some $C>0$} (see the definition 
and the details in \cite{2006}). The notion of beta-mixing -- which is neither defined nor discussed here  in detail -- is rather close  although not identical to 
the convergence in total variation. Hence, and also because of close results about lower bounds for 
convergence rates in \cite{Klokov_lb, MenshPopov}, it is likely that convergence of $Y$ to the stationary 
distribution of $Y_t$ under (\ref{eq:dens}) is also no faster than some  polynomial.  The assumption (\ref{binfty}) will not be used in the sequel but was shown just for information.  Recall that our aim is a faster convergence, and that we want to avoid any conditions on the derivative $\pi'$ except for its existence and local boundedness. \\ 

Note that in the case if for large \(|x|\) the density equals {\em exaclty} \(c(1+|x|)^{-m}\), it apparently follows that we need 
\(m>3\) in order to have the inequality  \(r>3/2\) in (\ref{binfty}). 
Yet, we will not use conditions in terms of \(\pi'\), assuming just (\ref{eq:dens}). Also, emphasize that the requirement \(m>3\) is for the quick reference on some existing earlier results. We do not claim that for \(m\le 3\) a similar analysis and asymptotics are not possible, but just that we are not aware of such asymptotics for \(m\le 3\).

Note that for the density on a half-line $\mathbb R^1_+$ a natural analogue of (\ref{eq:lang}) is the process satisfying an SDE with a non-sticky reflection at zero, 
\begin{equation}\label{eq:lang2}
dY_{t}=dW_{t}+b\left(Y_{t}\right)dt + d\phi^Y_t, \quad Y_0=\xi.
\end{equation}
For the process satisfying (\ref{eq:lang2}) similar mixing and convergence bounds follow from the bounds and from the calculus quite similar to those in \cite{ayv_grad_drift} applied to the situation of the reflected SDE, or just from a consideration of an SDE (\ref{eq:lang}) with a symmetric
(\(b(-x)=b(x)\)) drift. 

\subsection{Receipt 2: SDE with an affine ``mean-reverted'' drift}\label{sec22}
There is another receipt different from (\ref{eq:lang}) offered in \cite{Bibby2005} 
(see also the references therein concerning some other earlier constructions): in a slightly 
simplified form it suggests to consider an SDE on the line
\begin{equation}\label{sde_bibby}
dZ_t = -(Z_t-\mu) \,dt + \sqrt{v(Z_t)}\,dW_t, 
\end{equation}
with an appropriate initial distribution (e.g., stationary $\pi$ as in the reference paper), 
with 
\begin{equation}\label{bibby2}
v(z) = \pi(z)^{-1} \, \int_{-\infty}^z (\mu-s)\,\pi(s)\,ds, \quad 
\mu = \int s\,\pi(s)\,ds.  
\end{equation}
It is, of course, assumed that $\mu$ is finite, and then it is easily proved that $v\ge 0$, so that the SDE (\ref{sde_bibby}) is well-defined. (Indeed, $\int_{-\infty}^\mu (\mu-s)\,\pi(s)\,ds \ge 0$ since $\mu-s \ge 0$ for $s\le \mu$; and for $s>\mu$ the values $\mu-s$ are negative but the whole integral $\int_{-\infty}^{+\infty} (\mu-s)\,\pi(s)\,ds = 0$, so that for any $z<+\infty$, $v(z)\ge 0$ as required.) 
However, of course, ``good'' ergodic properties of the solution of this equation depend on 
some features of the density $\pi$. The solution locally exists due to local Lipschitz property 
of $\sqrt{v}$ combined with the affine drift assumption, but no-explosion should be derived from 
other conditions. Some related papers are, for example,  \cite{Abu-Ver09, Abu-Ver09b, kulik-leonenko} which 
tackle particular parametric families of target densities $\pi$ -- Student, reciprocal Gamma, and 
Fisher-Snedekor diffusions. 
In all three papers 
a quadratic Lyapunov 
function allows to show an exponential convergence in total variation which is non-uniform in the 
initial state $Z_0$. It is interesting that in \cite{Bibby2005} an 
exponential character of the stationary correlation function is established; yet, 
convergence of a {\em non-}stationary process to a stationary regime was not studied. 
In fact, the process under investigation in \cite{Bibby2005} is stated to be ``ergodic'', 
which ergodicity is understood in the sense of being stationary without any convergence 
statements. At the same time, the assumptions on the density in \cite{Bibby2005} involve the 
stationary density \(\pi\) (in our notations) but not on its derivative (possibly with 
some additional non-restrictive requirements in some theorems like continuity of the 
target density). Recall that in the present paper $C^2$-differentiability of \(\pi\) is assumed, 
but convergence rate bounds only depend on the asymptotic assumptions at infinity on the 
density \(\pi\) itself. It looks plausible that, in principle, it may be possible to work with ``weak'' definitions of 
the process via Dirichlet forms theory (\cite{Fuku, MaRo}), 
but we prefer to have a well-defined solution trajectory; in 
particular, we will be working with strong  solutions due to \cite{Ve81}. 
The receipt 2 naturally rises the question whether it is possible to arrange even a faster convergence, let theoretically.

\section{The setting \& main result}
Our primary goal is a density $\pi$ on $\mathbb R^1_+ = [0,\infty)$ satisfying
\begin{equation}
	c\left(1+x\right)^{-m}\leq
		\pi\left(x\right)\leq 
			c^{-1}\left(1+x\right)^{-m}, \quad x\ge 0.
\label{eq:dens1}
\end{equation}
Receipts I \& II in the previous section can be applied to this setting if we just extend the density in a symmetric way to the whole line, with a natural normalisation. As was said in the Introduction, we aim at constructing a diffusion process on $\mathbb R^1_+$ which converges towards $\pi$ with an exponential rate uniformly with respect to the initial data. It is likely that  a similar result holds true for a {\it symmetric} density $\pi$ on the whole  line $\mathbb R^1$ satisfying (\ref{eq:dens}), which we mention as a remark.

In order to achieve yet a  better convergence than typically guaranteed by the receipts (\ref{eq:lang})--(\ref{eqb}) or even by (\ref{sde_bibby})--(\ref{bibby2}), and for yet a more general class of densities than in \cite{Abu-Ver09, Bibby2005, kulik-leonenko} and in quite a few other works, let us consider two diffusion processes $Y_t$ and $X_{t}$ on $\mathbb R_+$ satisfying, respectively, SDEs with a non-sticky reflection at zero (\ref{eq:lang2})
and
\begin{equation}\label{feq}
dX_{t}=f\left(X_{t}\right)dW_{t}+f^{2}\left(X_{t}\right) b\left(X_{t}\right)dt + d\phi^X_t, \quad X_0 = \xi,
\end{equation}
with a local time $\phi_t^X$ at zero and 
with a special auxiliary function $f$, 
\begin{equation}
	f(z):= 
			\displaystyle \left(1+\int_{0}^{z}\frac{dy}{\pi(y)}\right)^{1/2},   \quad z\ge 0.
\label{eq:accel}\end{equation}
The generator  of this process is given by 
\[
 	L=f^2 L_0,
 \]
 where $L_0$ is the generator of the reflected  diffusion \eqref{eq:lang2}: 
\[
L_0v(x) = \frac12 v''(x) + b(x) v'(x), \; \forall \, x>0, \quad \& \quad L_0v(0) = v'(0+). 
\] 
Recall the requirements on the non-sticky solution and on its local time: $\phi^X$ is a monotonically non-decreasing function; for any $t>0$, 
\begin{align*}
\phi^X_t = \int_0^t 1(X_s=0)d\phi^X_s; \quad \int_0^t 1(X_s=0)ds = 0 \; \mbox{a.s.}
\end{align*}

~

Of course, a question about existence of solution of this equation (\ref{feq}) 
on the whole half-line \(t\ge 0\) arises here, and a positive answer to this question for the first sight may look doubtful given fast increasing coefficients. However, it will be 
justified with the help of a random time change and of the law of large numbers that such a (strong) solution exists on the whole line and does not explode. 
The main result  
is the following Theorem.

\begin{theorem}\label{thm1} 
Assume that for a strictly positive  probability density $\pi\in C^2$ with two locally bounded derivatives the bounds \eqref{eq:dens1}
hold with some $m>3$. 
Then the SDE (\ref{feq}) has a strong solution $X_{t}$ for all \(t\ge 0\) which is  
strongly (pathwise) unique and which possesses an exponential
rate of convergence to the stationary distribution $\pi\left(x\right)dx$, 
\begin{equation}\label{eq:crate}
\|\mu_t^{\xi} - \mu\| _{TV} \le C \exp(-\lambda t), \quad t\ge 0, 
\end{equation}
uniformly with respect to $\xi$, with some constants $\lambda$ and $C$ which both admit certain evaluation, where $\mu^\xi_t$ is a marginal measure of the process $X_t$ that starts from $\xi$ at $t=0$, and $\mu (dx) = \pi(x)dx$ is the (unique) invariant measure of the process. 

\end{theorem}
The right hand side in (\ref{eq:crate}) does not depend on the initial value $\xi$. 
Theoretical evaluations of both constants in the bound (\ref{eq:crate}) is likely to be not very efficient, yet possible which is clearly better than pure existence of such constants.

\section{Proof of Theorem \ref{thm1}}
The proof will be split into several steps.

\medskip 

\noindent
{\bf Step 1. Random change of time.}
Define the function $f(z)$ on 
 $\mathbb{R}^1_+$ by \eqref{eq:accel}.
Obviously there exists $0<a\le 1$ (namely, any $a\in [0,c^2]$, with $c$ from (\ref{eq:dens1})) such that 
\begin{equation}
	a\left(1+z\right)^{m+1}\leq f^{2}(z)\leq a^{-1}\left(1+z\right)^{m+1}, \quad \forall z \in \mathbb{R}_+^1.
\label{eq:est_f}\end{equation}
Let us define a random time change (cf. \cite{GS, McKean}) by 
\begin{equation}
	\chi_{t}:=\int_{0}^{t}f^{-2}\left(Y_{s}\right)ds, \quad \& \quad \beta_t:= \chi^{-1}_t \quad \mbox{(the inverse function).}
\label{eq:time-change}
\end{equation}
In other words, 
\[
\beta'_t = f^2(Y_{\beta_t}), 
\]
and
\[
t = \int_{0}^{\beta_t}f^{-2}\left(Y_{s}\right)ds.
\]
This time change $t\mapsto \beta_t$  is non-degenerate, that is, the following two conditions hold:

$\rm{(i)}$ there is no blow up at finite time: 
\begin{equation}
	\mathbb{P}(\chi_{t}|_{t\rightarrow T-0}\rightarrow+\infty) =0\quad \forall T\in(0,+\infty).
\label{eq:tc-bounds-2}\end{equation}

$\rm{(ii)}$ $\chi_t$ is unbounded as $t\rightarrow +\infty$ 
(i.e. when "real" time goes to infinity): 
\begin{equation}
	\chi_{t}\geq0,\quad\chi_{t}|_{t\rightarrow\infty}\rightarrow+\infty\quad\mathbb{P}-\mbox{a.s.}
\label{eq:tc-bounds-1}\end{equation}

To prove  \eqref{eq:tc-bounds-2}, it suffices to notice that, due
to \eqref{eq:est_f}, for any \(s<t\),
\[
	0\le \chi_{t}-\chi_{s}\leq a^{-1}\int_{s}^{t}\left(1+\left|Y_{r}\right|\right)^{-m-1}dr,
\]
hence 
$$
\chi_{t}^{'}\leq a^{-1}\cdot\sup_{r\in\mathbb{R}^{1}}\left(1+|r|\right)^{-m-1}=a^{-1}, \quad \mathbb{P}-\mbox{a.s.}
$$ 
Then (\ref{eq:tc-bounds-2}) immediately follows.

~

From here we find, 
\[
\inf_{t\ge 0}\beta'_t \ge a>0, \quad \mathbb{P}-\mbox{a.s.}, 
\]
and 
\[
\mathbb P(\limsup_{t\to\infty} \beta_t < \infty) = 0. 
\]

~

The  assertion \eqref{eq:tc-bounds-1} follows from the following Lemma. 
\begin{lemma}\label{erg} Let $m>3$, and let $g$ be a bounded continuous function
on $\mathbb{R}^{1}$. Assume that the diffusion process $Y_{t}$ satisfies
 (\ref{eq:lang}),  and let $\mu_{inv}$
be its unique invariant measure. 
Then for any $\delta>0$ and $\varepsilon>0$ there exists $T_0>0$
such that 
\[
	\mathbb{P}\left(\left|\frac{1}{t}\int_{0}^{t}g\left(Y_{s}\right)ds-\int g(x)d\mu_{inv}(x)\right|
		>\varepsilon\right)<\delta\qquad\mbox{for any }t\geq T_{0}.
\]
\end{lemma} 
The  Lemma with $g(r)=\left(1+\left|r\right|\right)^{-m-1}$ yields 
the  assertion (\ref{eq:tc-bounds-1}).
Indeed,
let us fix any $\delta\in(0,1)$ and \(\varepsilon = a_g/2\). Naturally, \(a_g > 0\). 
Then with $\mathbb{P}$-probability at least $1-\delta$
one has $\chi_{t}\geq\left(a_{g}-\varepsilon\right)t=a_g \, t/2$ for {\bf any} $t$ large
enough. This means that  with probability at least \(1-\delta\) the change of time mapping does not stop up to at least $a_g\, t/2$.  Since \(\delta\in (0,1)\) is arbitrary, \eqref{eq:tc-bounds-1} holds. 

\medskip

Proof of Lemma \ref{erg}. First of all, we will refer to the mixing results for SDEs on the whole line; however, in the case of symmetric coefficients ($b(-x)=b(x)$ and similarly for the diffusion if it is not a constant) such results straightforward imply similar bounds and convergence rates for (non-sticky) reflected at zero diffusions, too.
In other words, 
The process \((Y_t)\) is Markov ergodic with a finite variance (and, in fact, with any moment  $m'<m-2$; $g$ is bounded) with a polynomial beta-mixing rate as well as convergence in total variation \(\beta^\xi(t) + |\mu_t^\xi - \mu_\infty|_{TV} \le C_k(\xi)(1+t)^{-k}\) with some \(C(\xi)\) for any \(k<m-1\), to the stationary regime $\mu_\infty$, see \cite{ayv_grad_drift, ayv_grad_drift2}. Indeed, the assumptions of \cite{ayv_grad_drift} are met with \(p=m-1\) where \(p\) is the standing parameter in \cite{ayv_grad_drift}. The assumption $m>3$ implies $k>2$. Moreover, the function $g$ is bounded; hence, the process $\int_0^t g(Y_s)\,ds$ possesses all moments (including exponential with any constant, although, this is far too much for our goal). The beta-mixing coefficient dominates the alpha-mixing, while certain convergence rate to zero of the alpha coefficient is the standing assumption in the Theorem 18.5.4 of \cite{IbrLinnik}. 
Hence, for the stationary regime, the assertion
of the Lemma  -- LLN -- follows from the Central Limit Theorem  \cite[Theorem 18.5.4]{IbrLinnik}.  Indeed, splitting
the integral from zero to $t$ into a sum $\sum_1^{[t]}$ plus $\int_{[t]}^t$,  the claim follows.  For a {\it nonstationary}
regime the desired LLN follows again from the CLT for the {\it stationary} case, from the Markov property, 
and from the polynomial convergence of $\mbox{Law }(Y_{s})$ to $\mu_{inv}$
in total variation, similarly to the proof of `` non-stationary CLT'' in \cite[Theorem 4]{Ve2} with the help of the results from  \cite[Theorem 1]{ayv_grad_drift} with $m>3$. After mixing bounds have been found, see also \cite{Ve_LN} for LLN (formally, in \cite{Ve_LN} mixing is exponential, but obviously any polynomial would do such that the related sums or integrals converge). 
 This finishes the proof of the Lemma~\ref{erg}. \hfill{}$\square$

\medskip

See also \cite{Eva} for close results under slightly different assumptions. As was already mentioned, the statement of the Lemma will be used straight away for our reflected diffusion~(\ref{eq:lang2}). 

~

\noindent
{\bf Step 2. Constructing the process $X_t$.}

On the probability space $\left(\Omega,\mathcal{F}, (\mathcal{F}_{t}),\mathbb{P}\right)$ with a solution $Y_t$ to the equation (\ref{eq:lang2}), 
let us introduce stochastic processes 
\[
	X_{t}:=Y_{\beta_{t}}, \quad  \phi^X_t = \phi^Y_{\beta_{t}}.
\]
Then due to the time change \cite[Theorem 3.15.5]{GS} 
it follows that the process $X_{t}$ satisfies an SDE
\begin{equation}
	dX_{t}=f\left(X_{t}\right)\, d W_{t}+f^{2}\left(X_{t}\right) b\left(X_{t}\right)dt + d \phi^X_t, \quad X_0 = \xi,
\label{eq:lang-ac}
\end{equation}
with a new Wiener process \(\displaystyle \tilde W_{t} = \int_0^{\beta_t} f^{-1}(X_s) \, dW_s\), and with the local time at zero \( \phi^X_t\); recall that $f(0)=1$. Indeed, outside zero the ``main part'' here 
$$1(X_t>0)dX_{t}=1(X_t>0)\left[f\left(X_{t}\right)\, d\tilde W_{t}+f^{2}\left(X_{t}\right) b\left(X_{t}\right)dt\right]$$ follows straightforward from \cite[Theorem 3.15.5]{GS}, and  
$$
1(X_t=0)dX_{t}=1(X_t=0)d \phi^X_t
$$ 
is a direct consequence of the equation 
$$
1(Y_t=0)dY_{t}\,=\,1(Y_t=0)d\phi^Y_t.
$$ 
Also, we have, 
\begin{align*}
\int_0^t 1(X_s=0)ds  = 0, \;\; 
\int_0^t 1(X_s=0)d\phi^X_s  = \phi^X_t.
\end{align*}
Finally, 
\[
X_t - \xi - \int_0^t f\left(X_{s}\right)\, d W_{s}+\int_0^t f^{2}\left(X_{s}\right) b\left(X_{s}\right)ds -  \phi^X_t = 0, \quad {\mbox{a.s.}}
\]
Thus, $X$ is the solution of the equation (\ref{eq:lang-ac})  with a non-sticky reflection, as required. 

~

The equation (\ref{eq:lang-ac}) can be also derived from the time change for the SDE (\ref{eq:lang}) on the whole line with a symmetric drift and symmetrically extended $f$ after the application of It\^o--Tanaka's formula to the modulus, 
\begin{equation}
	d\bar X_{t}=\bar f\left(\bar X_{t}\right)\, d\tilde W_{t}+\bar f^{2}\left(\bar X_{t}\right) b\left(\bar X_{t}\right)dt, \quad \bar X_0 = \xi,
\label{eq:lang-ac2}
\end{equation}
with 
\[
\bar b(x) = \mbox{sign}(x)b(|x|), \;\; \bar f(x) = f(|x|), \quad \forall \, x\in \mathbb R^1.
\]

By construction, the processes $X_t$ and $ \phi^X_t$ are regular, i.e. are defined for all $t\geq 0$, and  adapted to the filtration \(\tilde {\cal F}_t \equiv {\cal F}_{\beta_t}\), see \cite{GS}.  Recall the well-known fact that the new filtration \({\cal F}_{\beta_t}\) is well-defined because of the fact that for any \(t\), the random variable \(\beta_t\) is a stopping time. 

\medskip

Emphasize that the process $X_t$ is well defined on the whole half-line $t\ge 0$, it does not explode, and it neither reaches infinity from zero, nor vice versa (zero from infinity) over a finite time, 
all of these because of the construction via the time change. 

\medskip

\noindent
{\bf Step 3.} The solution $X_t$ is strong. Indeed, it is well-defined on \(t\ge 0\), and the diffusion coefficient is locally continuously differentiable, and locally bounded,  and locally non-degenerate, while the drift coefficient is also locally bounded. Due to the results in \cite{Ve81}, 
this suffices for strong uniqueness via the stopping time arguments  with the help of the strong Markov property -- see \cite{Krylov_selection}.  This will be used in the sequel  in the coupling procedure (although, probably could be done with weak solutions, too). 

\medskip 

\noindent
{\bf Step 4. Stationary distribution for  $X_t$.}

Let us prove that the process \(X_t\) has a unique invariant distribution $\pi(x)dx$. 
The stationary distribution $\mu$ satisfies the stationary Kolmogorov equation $L^* (\mu)=0$ on $\mathbb R_+$ -- or, equivalently, \(L^*\pi = 0\) -- 
where 
\begin{equation}\label{L}
L=\frac{f^2}{2} D_x^2 +(f^2 b)D_x
\end{equation} 
is the generator of $X_t$ and 
${}^*$ is the adjoint with respect to the Lebesgue measure. 
First of all, the Kolmogorov equation $L^* (\mu)=0$ has at most one probability solution due to 
\cite[Example 4.1.1]{Bo}. Next, the measure $\mu(dx)=\pi(x) dx$ 
satisfies this equation. Indeed, since 
\begin{equation}
	\frac{1}{2}\pi^{'}-\left(b\pi\right)=0,
\label{eq:pl}\end{equation}
we have   
\begin{multline}
	\frac{1}{2}\left((f^{2})\pi\right)^{''}-\left((f^{2})b\pi\right)^{'}=
	\frac{1}{2}\left((f^{2})^{''}\pi+2(f^{2})^{'}\pi^{'}+(f^{2})\pi^{''}\right)-
	(f^{2})^{'}(b\pi)-(f^{2})(b\pi)^{'}=\\
	(f^{2})\left(\frac{1}{2}\pi^{''}-
	\left(b\pi\right)^{'}\right)+(f^{2})^{'}\pi^{'}+
	\frac{1}{2}(f^{2})^{''}\pi-(f^{2})^{'}(b\pi)\overset{\eqref{eq:pl}}{=}
	\frac{1}{2}(f^{2})^{''}\pi+\frac{1}{2}(f^{2})^{'}\pi^{'}=
	\frac{1}{2}\left((f^{2})^{'}\pi\right)^{'}.
\label{eq:check}\end{multline}
 But 
\begin{equation}\label{one}
	((f^{2})^{'}\pi)(x)= 
	\displaystyle		\left(1+\int_{0}^{x}\frac{dy}{\pi(y)}\right)^{'}\cdot\pi(x)=
	\frac{\pi(x)}{\pi(x)}=1,   \quad x\ge 0, 
\end{equation}
where at zero derivative is understood as right one. Hence the expression in the right hand side of \eqref{eq:check} equals zero, i.e., $\pi$ is a stationary measure for the new process \(X\). Note that the same calculus with $f$ replaced by $1$ shows the invariance of $\pi$ for $Y_t$ 

\medskip

As may be expected, (\ref{one}) implies  the equality 
\begin{equation}\label{e-inv}
	\mathbb{E}_\pi h(X_t) = \int h(y)\pi(y)\,dy, 
\end{equation}
for any \(t>0\) and \(h\in C_b(\mathbb{R}^1_+)\). 
By virtue of the Lebesgue dominated convergence 
theorem we can take \(h\in C^\infty_0(\mathbb{R}^1_+)\) (continuous with a compact support), but more than that, it suffices to consider functions $h\in C^\infty_0(\mathbb{R}^1_+)$ with $h'(0+)=0$: indeed, the latter class -- denoted in the sequel as $C^\infty_{00}(\mathbb{R}^1_+)$ -- is clearly dense in $C^\infty_0(\mathbb{R}^1_+)$.  
We have, due to $h'(0+)=0$, 
\begin{align*}
	dh(X_t) = Lh(X_t)\,dt + h'(X_t)f(X_t)\,dW_t + 1(X_t=0)h'(0+)d\phi^X_t
	\\\\
= Lh(X_t)\,dt + h'(X_t)f(X_t)\,dW_t.	 
\end{align*}
Moreover, since \(h\) has a compact support,  \(Lh\) and \(h'f\) are bounded. So, by rewriting
in integral form and taking expectations we get, 
\[
	\mathbb{E}_\pi h(X_t)  -  \mathbb{E}_\pi h(X_0) = \mathbb{E}_\pi \int_0^t Lh(X_s)\,ds 
\equiv  \int_0^t \mathbb{E}_\pi Lh(X_s)\,ds, 
\]
the last equality due to Fubini's theorem. 
Denote by $p_s(y,z)$ the transition density of the Markov process $X$; its existence follows, e.g., from \cite[Corollary 2.9 \& Remark 2.17]{Bogach}. In a ``good case'' with all appropriate derivatives, the standard {\it formal} calculus runs as follows: 
\begin{align*}
	 	 \mathbb{E}_\pi Lh(X_s) = \iint (L_zh(z))\pi(y)p_s(y,z)\,dzdy
	\\\\
= \iint \pi(y) h(z) L^*_zp_s(y,z)dzdy
= \iint \pi(y) h(z) \partial_s p_s(y,z)dzdy
 \\\\
= -\iint \pi(y) h(z) L_y p_s(y,z)dzdy
= -\int h(z) \left(\int p_s(y,z) L^*_y \pi(y)dy\right)dz = 0,  
\end{align*}
due to 
forward and backward Kolmogorov's equations, and Fubini's theorem. Therefore, we conclude that 
\begin{equation}\label{esta}
\mathbb{E}_\pi h(X_t)  =  \mathbb{E}_\pi h(X_0), 
\end{equation}
which is equivalent to (\ref{e-inv}). 
A rigorous justification without additional assumptions follows from \cite{Bogach}.

Note that uniqueness of the invariant measure will follow from the convergence bound (\ref {eq:crate}) once it is established.

\medskip

\noindent
{\bf Step 5. Uniform exponential moment bound.}

\noindent
Let us take {\it any} $K>0$, and  define 
\begin{equation}	\gamma_{X}^{\xi}\equiv\gamma_X^{\xi,K} \equiv\gamma :=\inf\left(t\geq0:\,\, X_{t} 
	\leq K,
				\,\, X_{0}=\xi \ge 0\right).
\label{eq:moments}
\end{equation}
Let us show that  $\E\exp(\alpha\gamma_{X}^{\xi}) <+\infty$
 for $\alpha>0$ small enough, uniformly with respect to the initial state of the process. 
Denote $v_q \left(\xi\right):=
\mathbb{E}_{\xi}\gamma^{q}$, 
with a convention $v_{0}\equiv1$. 
 Obviously, 
\[
	\mathbb{E}_{\xi}e^{\alpha\gamma}	=\sum_{q=0}^{+\infty}\frac{\alpha^{q}\mathbb{E}_{\xi}\gamma^{q}}{q!} = \sum_{q=0}^{+\infty}\frac{\alpha^{q}v_{q}\left(\xi\right)}{q!} < \infty,
\]
provided all quantities \(v_{q}(\xi)\) are finite and grow not too fast  in $q$. 
 
\medskip

\noindent
{\bf Step 6.  Auxiliary results for polynomial moments.}

\noindent
In order to guarantee that the values $v_q$, indeed, may not grow too fast, let us find alternative representations for them. By
virtue of the identity 
\[	\left(\int_{0}^{\gamma}1dt\right)^{q}=q\int_{0}^{\gamma}\left(\int_{t}^{\gamma}1ds\right)^{q-1}dt,
\]
which holds both for finite and infinite $\gamma$, we get,
\[
v_{q}(\xi)=q\mathbb{E}_{\xi}\int_{0}^{\gamma}v_{q-1}(X_{t})dt
\]
at least, for any $q\ge 1$ such that $v_{q}\left(\xi\right)$ is finite. Indeed, due
to the Fubini's theorem (iii) and the Markov property (iv), 
\begin{multline*}	v_{q}\left(\xi\right)\equiv \mathbb{E}_{\xi}\gamma^{q} =	q \mathbb{E}_{\xi}\int_{0}^{\gamma}
\left(\int_{t}^{\gamma}1ds \right)^{q-1}dt =	q\mathbb{E}_{\xi}\int_{0}^{\infty}1 \left(\gamma > t\right)		\left(\int_{t}^{\gamma}1ds\right)^{q-1}dt =
 \\
\overset{(iii)}{=}q\int_{0}^{\infty}\mathbb{E}_{\xi}
		1 \left(\gamma >  t\right)\left(\int_{t}^{\gamma} 1ds\right)^{q-1}dt =
	q\int_{0}^{\infty}\mathbb{E}_{\xi}1 \left(\gamma\geq t\right)	\mathbb{E}_{\xi}\left(\left(\int_{t}^{\gamma}
	1ds\right)^{q-1} | {\cal F}^X_t\right) dt=
 \\
=
	q\int_{0}^{\infty}\mathbb{E}_{\xi}1 \left(\gamma >  t\right)	\mathbb{E}_{\xi}\left((\gamma - t)^{q-1} | X_t\right) dt	=q\int_{0}^{\infty}\mathbb{E}_{\xi}1 \left(\gamma >  t\right)
		\mathbb{E}_{X_{t}}\gamma^{q-1}dt 
=
 \\
\overset{(iv)}{=}q\int_{0}^{\infty}\mathbb{E}_{\xi}1 \left(\gamma >  t\right)v_{q-1}(X_{t})dt\overset{(iii)}{=}	q\mathbb{E}_{\xi}\int_{0}^{\gamma}v_{q-1}(X_{t})dt.
\end{multline*}
Hence, if both quantities are finite, then 
\[	v_{q}\left(\xi\right)=q\mathbb{E}_{\xi}\int_{0}^{\gamma}v_{q-1}(X_{t})dt,\quad v_0=1.
\]
Note that for each $q\geq1$, if $v_q$ is finite, then it satisfies  
\begin{equation}
	Lv_{q}(x)=-qv_{q-1}(x), \;\; x\ge K, 
\label{eq:pois}\end{equation}
by virtue of the probabilistic representation of solutions of the elliptic equation with Dirichlet boundary condition,  or equivalently by Duhamel's formula. 
Obviously  $v_{q}(K)=0$, as well as $v_{q}(x)=0, \, 0\le x\le K$. Also, 
it is known that  if $v_q(\xi)<\infty$ for some $\xi$ then it is finite for any $\xi$. 
However, we 
are not going to use this equation directly since it lacks the ``second boundary condition'' normally required for the second order differential equation. Instead, we will find solutions to boundary problems that  approximate $v_q$. 
In fact, what we shall need instead is the following Lemma.

Let  $\hat L:\,(\hat Lu)(x)=a(x)u^{''}(x)+c(x)u^{'}(x)$ be the generator of the diffusion process $(\zeta_t,t\geq 0)$ 
with locally bounded coefficients $a>0$ (the diffusion) and $c$ (the drift), and such that $a$ is locally uniformly non-degenerate, and which process is a strong solution of a corresponding SDE.
Let us fix a positive (non-negative) function $\psi$ on $\mathbb{R}^1$. Let 
\[
	\tau_K:=\inf\left(t\geq 0:\,\,\zeta_{t}\leq K,
				\,\, \zeta_{0}=\xi\right) \quad (\xi>0),
\]
and 
\[	
v\left(\xi\right) =\mathbb{E}_{\xi}\int_{0}^{\tau_K} \psi(\zeta_{t})\,dt.
\]
\begin{lemma}\label{PDE_lemma} 
For any 
\(N>K>0\), consider the boundary problem 
\begin{equation}
	\hat Lv^+_{N}=-\psi,\quad v^+_{N}\left(K\right)=0,\,\,\left(v^+_{N}\right)^{'}(N)=0,
\label{eq:bvp-lem}
\end{equation}
Then the function $v^+_{N}(\xi)\uparrow v(\xi)$  as $N\uparrow\infty$, for every
$\zeta_0 = \xi$ with $\xi \ge K$.
\end{lemma} 
 \textbf{Proof of Lemma \ref{PDE_lemma}.} 
For any 
 $0\le K\le \xi\le N$, let us consider a family of stochastic
processes $\zeta_{t}^{N}$, given by the SDE with reflection, 
\[
	d\zeta_{t}^{N}=\sqrt{2a(\zeta_t ^N)}dw_{t}+c(\zeta_{t}^{N})dt +d\phi_{t}^{N}, \quad \zeta_0^N = \xi, 
\]
with values on $[0,N]$, 
with a non-sticky reflection at $N$ and an absorbtion at zero, where $\phi_{t}^{N}$ is its local time at $N$.
Applying It\^o's formula (or, in fact, more precisely It\^o--Krylov's formula if  continuity of $\psi$, $a$, and $c$ is not assumed) to $v_{N}(\zeta_{t}^{N})$,
we get the following representations:
\[
	v_{N}(x)=\E_{x}\int_{0}^{\tau_{K,N}}\psi(\zeta_{s}^{N})ds,
\]
where 
\(
	\tau_{K,N} = \inf\left(t\geq0:\,\, |\zeta^N _{t}| \le K\right)
\) is the moment when the process $\zeta^N _t$ first hits the interval $[0,K]$. 
Note that 
$\tau_{K,N}$ monotonically increases as $N$ increases. Also note that  if $\tau_K<\infty$ then, obviously, \(\tau_{K,N} \uparrow \tau_K\); and if 
$\tau_K=\infty$ then still \(\tau_{K,N} \uparrow \infty =\tau_K\). These all follow   from the comparison theorem for one-dimensional SDEs possessing strong solutions with the same coefficients and different initial data. This comparison theorem can be shown as follows. Consider two SDEs with the same initial value $\xi$ but with two different $N_1 < N_2$, say. Denote the corresponding solutions by $\zeta^{N_1}_t$ and $\zeta^{N_2}_t$. Assuming that $\xi\in [K,N_1]$, due to the strong uniqueness $\zeta^{N_1}_t = \zeta^{N_2}_t$ until $\hat \tau^{}:=\inf(t\ge 0: \, \zeta^{N_1}_t = K \; \mbox{or} \; N_1)$. If at this moment -- which is a stopping time -- 
$\zeta^{N_1}_t = \zeta^{N_2}_t = K$, then the claim is justified because $K$ is the absorbtion point. If, however, $\zeta^{N_1}_t = \zeta^{N_2}_t = N_1$, then the first process $\zeta^{N_1}_t$  will remain less than or equal to $N_1$ all the time, while the second will exceed this level $N_1$ with probability one on any right interval of $\hat \tau$. This follows easily from the ``reverse'' time change which makes diffusion back equal to one and from the Girsanov theorem about eliminating the drift via change of measure, because for the standard Wiener process this property is well-known (e.g., it follows from Khintchin's iterated logarithm law for WP \cite{ItMc}
along with the strong Markov property. Thus, on any small right neighbourhood of the moment $\hat\tau$ we would have $\zeta^{N_1} \le \zeta^{N_2}$, with strict inequality at least at infinitely many  moments of time arbitrarily close to $\hat \tau$. Yet, both solutions are strong Markov. So, if we now start two processes with the same generator a new at two distinct initial value $\xi_1 < \xi_2$, then due to continuity the two solutions will satisfy 
$1(t>\hat\tau)1(\zeta^{N_1}_t<\zeta^{N_2}_t) = 1(t>\hat\tau)$, at least, until they meet again, i.e., for all $t< \bar\tau:= \inf(s\ge \hat\tau: \zeta^{N_1}_s = \zeta^{N_2}_s)$ (here, of course, $\inf(\emptyset) = \infty$, and, in fact, they will never meet again). But then, if we assume that $\bar \tau < \infty$, they will again coincide until the next moment when they touch the level $N_1$, after which we have again $\zeta^{N_1}\le \zeta^{N_2}$, and the cycle can repeated indefinitely. This shows that  
$\zeta^{N_1}_t\le \zeta^{N_2}_t$ for all $t\ge 0$. 
Hence, we have \(\tau_{K,N} \uparrow \tau_K\), $N\uparrow \infty$, and so, the monotonic convergence Theorem yields  the assertion of the Lemma, as required. 
 \hfill{}$\square$ 

\noindent
Similar calculi in similar situations yielding various close claims can be found in \cite{Mao, ayv_grad_drift2, ayv_grad_drift}.

~

\noindent
{\bf Step 7. Bounds for polynomial moments.}

\noindent
Let us prove that 
\begin{equation}
	v_{q}\left(\xi\right)\leq q!\cdot C^{q}
\label{eq:major_est-1}
\end{equation}
for all $q\geq1$, with 
\[
	C=\frac{a^{-1}}{m}\cdot A_{m},\quad A_{m}:=\int_{K}^{\infty}\left(1+w\right)^{-m}dw.
\]
does not depend on $q$. 
We argue by induction. Now  our particular generator is $\hat L = L$ from~(\ref{L}). 

\medskip 

\underbar{Base}: Let $q=1$. Fix $N>0$. Notice that $v^{0}=1$ and
consider a boundary value problem,
\begin{equation}
	Lv_{N}^{1}=-1,\quad v_{N}^{1}\left(K\right)=0,\,\,\left(v_{N}^{1}\right)^{'}(N)=0.
\label{eq:BVP-1}
\end{equation}
Since $L=f^{2}L_{0}$,
where $L_{0}u=\frac{1}{2}u^{''}+\frac{1}{2}\nabla\ln\pi(x)u^{'}$, this problem
admits a unique solution 
\begin{equation}
	v_{N}^{1}(\xi)=2\int_{K}^{\xi}\pi^{-1}(w_{1})dw_{1}
		\int_{w_{1}}^{N}\frac{\pi(w_{2})}{f^{2}\left(w_{2}\right)}dw_{2},
\label{eq:sol_pois__no_lim-1}\end{equation}
and due to \eqref{eq:dens} and \eqref{eq:est_f} we estimate replacing \(N\) by infinity in the upper limit of the integral, 
\begin{multline*}
	v_{N}^{1}(\xi)\leq2c\int_{K}^{\xi}(1+w_{1})^{m}dw_{1}\int_{w_{1}}^{\infty}c^{-1}\cdot
		(1+w_{2})^{-m}\cdot a^{-1}\left(1+w_{2}\right)^{-m-1}dw_{2}\leq\\
	\leq2a^{-1}\int_{K}^{\xi}(1+w_{1})^{m}dw_{1}\int_{w_{1}}^{\infty}(1+w_{2})^{-m}\left(1+w_{2}\right)^{-m-1}dw_{2}
	=\\
	=2a^{-1}\int_{K}^{\xi}(1+w_{1})^{m}\frac{(1+w_{1})^{-2m}}{2m}dw_{1}\leq
		\frac{a^{-1}}{m}\int_{K}^{\infty}(1+w_{1})^{-m}dw_{1}=A_{m}\cdot\frac{a^{-1}}{m}=:C.
\end{multline*}
By virtue of Lemma \ref{PDE_lemma},  this implies 
$v_{1}\left(\xi\right)=\lim_{N\rightarrow\infty}v_{N}^{1}\left(\xi\right)\leq C$.

\medskip

\underbar{Induction Step}: Note that if the right hand side in the equation (\ref{eq:BVP-1}) 
is multiplied by a constant \(R>0\), then, given the specific boundary conditions,  the bound for the solution will be also multiplied by this \(R\), so that instead of the upper bound \(C\) there will be a new upper bound \(RC\). 

Suppose that for some \(q\) and  for $n=q-1$ we have 
\[
	v_{n}(\xi)\leq C^{n}\, n!
\]
with the same constant \(C\) as above. Then, by the remark in the beginning of the induction step with \(R=C^{n}\, n! \times q \equiv C^{q-1}\,q!\), we immediately obtain
\[
	v_{q}(\xi)\leq C^{q-1}\, q!\times C = C^{q}\, q!, 
\]
as required. Hence, the inequality (\ref{eq:major_est-1}) follows. 
Note that a similar simple argument with a reference to the induction method and without using explicitly the second barrier \(N\) can be found in \cite[Lemma 3.1]{Mao}; practically the same calculus, yet with unbounded  growing in \(x\) moments was used in \cite{ayv_grad_drift}.

~
 
Now, take any $\alpha\in\left(0,C^{-1}\right)$. Then due to (\ref{eq:major_est-1})
one has 
\begin{equation}
\mathbb{E}_{\xi}e^{\alpha\gamma}=\sum_{q=0}^{+\infty}\frac{\alpha^{q}\mathbb{E}_{\xi}\gamma^{q}}{q!}
		\leq\sum_{q=0}^{+\infty}\alpha^{q}
		C^{q}=\frac{1}{1-\alpha C}<\infty.
\label{eq:exp}
\end{equation}
It may be argued now that the desired ``exponential coupling'' can  be arranged via the exponential moment bound (\ref{eq:exp})  and due to the elliptic Harnack inequality for divergent type equations \cite[Theorem 8.20]{GT} in the way similar to \cite{ayv_grad_drift}, see the next step. Note that it is a ``common knowledge'' that the bound (\ref{eq:major_est-1}) suffices for the Theorem claim. The reader who knows the exact reference may skip the rest of the proof.

\medskip 

\noindent
{\bf Step 8. Using Harnack inequality. }  
The usage of coupling method assumes that glueing or meeting of two versions of the process -- one stationary and another non-stationary -- can be arranged with a positive probability bounded away from zero on each  period of this construction. Here it suffices to consider a ``symmetric'' SDE on the whole line. By ``period'' in our case any finite interval may be taken; e.g., it is convenient to use \([0,2]\) which will be split into two equal parts, \([0,1]\) and \([1,2]\) (their intersection at one single point is not important). On the first half, according to the inequality (\ref{eq:exp}) and Bienaym\'e -- Chebyshev -- Markov's inequality, both independent versions of the process will attain some (actually, any) bounded neighbourhood of zero. On the second half we want to glue them with a probability also bounded away from zero. Note that the standard 1D or finite state space idea just to wait until the two trajectories intersect here does not work straightforward as we would like it. Or, rather, it works but the bound obtained in such a way would use some bounds on the derivative \(\pi'\), which we want to avoid by all means. 

There is a recent rather general tool based on regeneration period moments \cite{Zv}. Yet, to verify the mild condition (*) required for this tool is probably no easier than -- or, maybe, equivalent to -- what we suggest instead in the next paragraphs. 

One more approach which does not involve any properties of \(\pi'\)  uses classical inequalities for divergence form PDEs.  
Indeed, this step justifies that it is possible by virtue of Moser's Harnack inequality for divergent type
elliptic equation (see \cite{GT}) (cf., e.g., \cite{Ve1}, \cite{Ve2} where a parabolic Harnack inequality was used for the same goal). 
Here it is convenient to return to an SDE on the whole line with symmetric coefficients (\ref{eq:lang-ac2}) which solution is denoted by $\bar X_t$. Its modulus satisfies the equation (\ref{eq:lang2}) with a new Wiener process. 

We argue that an elliptic Harnack inequality
\begin{equation}\label{ha}
\mathbb E_x g(\bar X_{\bar\sigma}) \le C \mathbb E_{x'} g(\bar X_{\bar\sigma})
\end{equation}
for any non-negative function $g$ and any $|x|, |x'| \le 1$ with
$\bar \sigma := \inf\,(t\ge 0:\; |\bar X_t|\ge 2)$ follows from  \cite[Theorem 8.20]{GT} due to the equation $\mbox{div} (\exp(2U(x) \nabla v(x)) = 0$, 
here $v(x) = \mathbb E_x g(\bar X_\sigma)$. This reasoning should be combined with the bound $\mathbb P_x(\bar\sigma>t) \le 
C t^{-1}$ for \(t>0\) with some $C>0$ depending on the sup-norms of all coefficients in the ball \(B:= \{|x|\le 2\}\). The latter bound follows from \cite[Theorem 8.16]{GT} applied and from the Bienaym\'e -- Chebyshev--Markov inequality \(\mathbb P_x(\sigma>t) \le t^{-1}\,\mathbb E_x \sigma\) since the function \(v(x) := \mathbb E_x \bar \sigma\) is a solution to the equation \(\frac12\,\exp(-2U)\mbox{div}(\exp(2U) \nabla v) + 1 = 0, \; \& \; v|_{\partial B}=0\), or, equivalently, to the (Poisson) equation
\[
\mbox{div}(\exp(2U) \nabla v) + 2\exp(2U) = 0, \quad v|_{\partial B}=0, 
\]
to which the Theorem 8.16 \cite{GT} is applicable stating that solution \(v(x)\) is bounded by a constant, say, \(N\) depending only on \(\sup_{|x|\le 2}|U(x)|\) (actually, even on some integral norm of \(\exp(2U)\)). This immediately implies that by choosing \(t\ge 3N\) we have that \(\mathbb P_x(\bar\sigma\le t) \ge \frac23\). The same estimate holds true for the process $X_t$, with the stopping time $\sigma = \inf\,(t\ge 0:\; X_t\ge 2)$, and with a non-negative function $g$ on $\mathbb R^1_+$, i.e., 
$$
\mathbb E_x g(X_{\sigma}) \le C \mathbb E_{x'} g(X_{\sigma}).
$$

Along with 
(\ref{eq:exp}), this suffices for a successful exponential coupling for the process $X_t$ with its stationary version.  Although it will not be used here, note that the obtained bound implies a stronger exponential inequality $\mathbb P_x(\sigma>t) \le 
C \exp(-\lambda t)$ with some $C,\lambda>0$ by the well-known property of homogeneous Markov processes and their exit times.  

Finally,  
we can change the function
$U$ outside the ball $|x|\le 3$ so that it becomes bounded, --
the latter is possible without changing the process until
$\sigma$.

\medskip

\noindent
{\bf Step 9. Exponential convergence.}  Let us return to the half-line $\mathbb R_+$ and to the process $X_t$, and let us fix some $K>0$. It is known that -- modulo the conclusion of the previous step -- for the proof of  the desired exponential convergence in total variation, it sufficies to show that $\E\exp(\alpha\gamma_{X}^{\xi})$ is finite
for some $\alpha>0$, 
and for some -- {\it actually, for any} -- $K>0$ where $\gamma_{X}^{\xi}= \gamma_{X}^{\xi, K}$ was defined in (\ref{eq:moments}).

Let $X_{t}^{st}$ be the {\it independent}  stationary version of the Markov process $X_{t}$, i.e., $X_{t}^{st}$  is the process with the same generator and initial distribution with the density \(\pi\), if necessary, on some extended probability space with another independent Wiener process. (However, we will not change our notations for $\mathbb P$ and $\mathbb E$.) Naturally, the couple $(X_{t}, X_{t}^{st})$ is considered on some extension of the original probability space.  For $\xi > K$ let $\tau\equiv \tau^\xi$ be the 
moment of the first intersection of $X_t$ started from \(X_0=\xi\) with the  stationary version $X_t ^{st}$, i.e., 
$$
	\tau :=\inf\left(t\geq0:\,\, X_{t}=X_{t}^{st}\right).
$$ 
 As the random variable
$\tau$ is a stopping time and $X_{t}$ has strong Markov property,
we can define a new strong Markov process 
\begin{equation}
\hat{X}_{t}:=X_{t}1\left(t<\tau\right) +X_{t}^{st}1\left(t\geq\tau\right),
\label{eq:coupling}\end{equation}
with the property
\[
\mbox{Law }(\hat X_{t}^{})=\mbox{Law }\left(X_{t}\right).
\]
Obviously on $\{t>\tau\} \cap \{\tau < \infty\}$ the trajectories of $X_{t}^{st}$ and $\hat{X}_{t}$ ``after $\tau$'' 
coincide. Then for any Borel set $A$ one has (we drop the initial value $\xi$ since the final estimate is uniform in it)
\begin{align*}
	\left|\mathbb P\left(X_{t}\in A\right)-\mathbb P\left(X_{t}^{st}\in A\right)\right|\overset{\eqref{eq:coupling}}{=}
	\left|\mathbb P\left(\hat{X}_{t}\in A\right)-\mathbb P\left(X_{t}^{st}\in A\right)\right|=
 \\\\
=\left|\mathbb{E}\left(\left(1 \left(\hat X_{t}\in A\right)-1 \left(X_{t}^{st}\in A\right)\right)
\times\left(1 \left(t<\tau\right)+1 \left(t\geq\tau\right)\right)\right)\right|=
 \\\\
=\left|\mathbb{E}\left(\left(1 \left(\hat X_{t}\in A\right)-1 \left(X_{t}^{st}\in A\right)\right)
	\times1 \left(t<\tau\right)\right)\right|
 \\\\
\leq\mathbb{E}\left|\left
	(1 \left(\hat X_{t}\in A\right)-1 \left(X_{t}^{st}\in A\right)\right)\right|
	\times1 \left(t<\tau\right).
\end{align*}
 Since $\left|\left(1 \left(\hat X_{t}\in A\right)-1 \left(X_{t}^{st}\in A\right)\right)\right|\leq1$,
taking into account the exponential version of Bienaym\'e -- Chebyshev--Markov's  inequality, we conclude that 
\[
	\left|\mathbb{P}\left(X_{t}\in A\right)-\mathbb{P}\left(X_{t}^{st}\in A\right)\right|\leq	\mathbb{E}1 \left(t<\tau\right)=\P\left(t<\tau\right)
	\leq\exp\left(-\alpha t\right)\E\exp(\alpha\tau).
\]
Passing to the supremum over Borel sets $A$ we obtain due to (\ref{eq:exp}),  
\[
	\left\Vert \mu_{t}-\mu_{st}\right\Vert _{TV}\leq2\exp\left(-\alpha t\right)\E\exp(\alpha\tau) \le \frac{2}{1-\alpha C}\,\exp\left(-\alpha t\right),
\]
where $\mu_{t}=Law\left(X_{t}\right)$ and $\mu_{st}(dx)=\pi(x)dx$. This bound does not depend on the initial value $\xi$ which was dropped in the notation $\mathbb E_\xi$. Hence, the proof of the Theorem is completed. 
$\hfill \square$

~

\begin{remark}
It is likely that for the symmetric density $\pi$ on $\mathbb R^1$ and for the equation (\ref{eq:lang}) on $\mathbb R^1$ this method is applicable, too, with the function $f(y) \equiv f(|y|)$, and that it provides a similar convergence rate bound (\ref{eq:crate}) as in the Theorem \ref{thm1}. We leave it till further papers.   
\end{remark}

~

\begin{remark}
The form of the process $X_t$ is a result of an educated guess. We were looking for the process $X_t$
 whose generator $L$ would be the 
generator of the Langevin diffusion multiplied by a positive function $F=f^2$, which  is needed for applying a random time change. At the same time, we wanted the 
new process $X_t$ to have the same  invariant density  $\pi$. 
From this condition the function $F$ is 
determined up to two constants $c_1, c_2>0$
\[
	F(x)=
	\displaystyle	c_1 + \int_0 ^{x} \frac{c_2 dy}{\pi(y)}, \quad  x\geq 0.
\]
It can be checked by an explicit computation that with this choice of the function \(F\), the new process 
$X_t$ would still have an exponential rate of convergence to the invariant measure (we choose $c_1=c_2=1$ but in fact any 
strictly positive $c_1$ and $c_2$ give the same result). 
\end{remark}

\begin{remark}
Note that in  the Theorem  \ref{thm1} the property of continuity of the state space is important.
 If the state space is discrete, the modification we consider (multiplication of the generator by a function) 
 typically does not affect the rate of convergence. Indeed, let us consider a birth-death process $Y_t,\,\,t\geq 0$
  with 
 birth rates $\{\lambda_n,\, n\in \mathbb{N}_{\geq 0}\}$ and death rates 
 $\{\mu_n,\, n\in \mathbb{N}\}$:
 \begin{equation}
	\mathbb{P}\left(Y_{t+h}=m|Y_{t}=n\right)=\left\{ \begin{array}{lcl}
		\lambda_{n}h+o(h), &  & m=n+1\\
		\mu_{n}h+o(h), &  & m=n-1\\
		1-\mu_{n}h-\lambda_{n}h+o(h), &  & m=n\\
		o(h), &  & |m-n|>1\end{array}\right.
\label{eq:BDP}\end{equation}
The generator $A_0$ of the process $Y_t$ is given by 
\[
	A_0 \varphi (n) = \lambda_n ( \varphi (n+1) - \varphi (n) ) + \mu_n (\varphi (n-1)-\varphi (n)).
\]
If we want to apply the same transformation as in the continuous case, i.e. to consider a 
birth-death process $X_t$ with the generator $A=f_n \cdot A_0$, then the new process $X_t$ will have  the 
birth and death rates  $\{\lambda_n ^{'}=f_n \cdot \lambda_n,\, n\in \mathbb{N}_{\geq 0}\}$ and 
$\{\mu_n ^{'}=f_n \cdot \mu_n,\, n\in \mathbb{N}\}$ respectively. 

The invariant distribution $\pi$ of $Y_t$ for the  can be computed
explicitly and equals
\[
	\pi(\{ n\} )\equiv\pi_{n}=
		\pi_{0}\frac{\lambda_{0}\dots\lambda_{n-1}}{\mu_{1}\dots\mu_{n}},
		\quad\pi_{0}=\left(1+\sum_{n\geq1}
			\frac{\lambda_{0}\dots\lambda_{n-1}}{\mu_{1}\dots\mu_{n}}\right)^{-1}.
\]
Hence 
\[
	 \frac{\lambda_{n-1}}{\mu_n} = \frac{\pi_{n-1}}{\pi_n}\quad \mbox{for each}\quad n\in\mathbb{N}.
\]
The assumption that $X_t$ has the same invariant distribution as $Y_t$ yields 
\[
	\frac{\lambda_{n-1} ^{'}}{\mu_n ^{'}}=\frac{f_{n-1} \cdot \lambda_{n-1}}{f_n \cdot \mu_n}  
		= \frac{\pi_{n-1}}{\pi_n}\quad \mbox{for each}\quad n\in\mathbb{N},
\]
hence $f_n = f_0 = \mbox{const}$ for all values of $n$. So, such a transformation is just changing the time scale by multiplying it by a positive constant which doesn't influence the rate of convergence  qualitatively. 

\end{remark}

\subsubsection*{Acknowledgements}
The authors are grateful to the anonymous referee for valuable remarks.

\bibliographystyle{plain}
\bibliography{MZ_base}

\begin{thebibliography}{99}

\bibitem{Abu-Ver09}
N.Abourashchi, A.Yu. Veretennikov,  On exponential mixing bounds and convergence rate for reciprocal Gamma diffusion processes, Mathematical Communications, 14(2) (2009), 331-339. 

\bibitem{Abu-Ver09b}
N. Abourashchi, A.Yu. Veretennikov, On exponential mixing and rate of convergence for Student processes, Theory Probab. Math. Statist. 81 (2010), 1-13.



\bibitem{AitSahalia}
Y. A\"it-Sahalia,  
Nonparametric pricing of interest rate derivative securities, Econometrica, 1996, 64(3), 527--560. 


\bibitem{Bibby2005}
B.M. Bibby, I.M. Skovgaard, M. S{\o}rensen, 
Diffusion-type models with given marginal distribution and autocorrelation function, 
Bernoulli, 11(2) (2005), 191--220.

\bibitem{Bogach}
V.I. Bogachev, N.V. Krylov, and M. R\"ockner,  On regularity of transition probabilities and invariant measures of singular diffusions under minimal conditions, Communications in Partial Differential
Equations, 26(11) (2001), 2037--2080.

\bibitem{Bo}
V.I. Bogachev, N.V. Krylov, M. R\"ockner, S.V. Shaposhnikov, 
Fokker--Planck--Kolmogorov Equations, AMS, Providence, RI, 2015.


\bibitem{Cattiaux}
P. Cattiaux,
A. Guillin, 
Hitting times, functional inequalities, Lyapunov conditions
and uniform ergodicity, 
Journal of Functional Analysis, 272 (2017) 2361-2391.

\bibitem{MP}
R. Douc, E. Moulines, P. Priouret, P. Soulier, 
Markov Chains, 
Springer, 2018.

\bibitem{Fort}
G. Fort, E. Moulines, P. Priouret, 
Convergence of adaptive and interacting Markov Chain Monte Carlo algorithms, 
Annals of Statistics, 39(6) (2011) 3262--3289. 


\bibitem{Fuku}
M. Fukushima, 
Dirichlet forms and Markov processes, 
Amsterdam; New York: North-Holland Pub. Co. \& 
Elsevier North-Holland, 1980.

\bibitem{GS}
I. Gikhman and A. V. Skorokhod, Stochastic Differential Equations, 
Springer, New York, 1972.

\bibitem{GT} 
D. Gilbarg, N.S. Trudinger,  Elliptic partial differential
equations of second order. 2nd ed.,  Springer, Berlin et al., 1983. 

\bibitem{IbrLinnik}
I.A Ibragimov, Yu.V. Linnik, 
Independent and stationary sequences of random variables, 
Groningen, Wolters-Noordhoff, 1971.


\bibitem{ItMc}
K. It\^o, H.P. McKean, 
Diffusion Processes and Their Sample Paths. 
Springer, 1965.

\bibitem{Klokov_lb}
S.A. Klokov, 
Lower Bounds of Mixing Rate for a Class of Markov Processes, 
Theory Probab. Appl., 51(3) (2007) 528--535. 

\bibitem{Khasminsky}
R. Khasminskii, 
Stochastic Stability of Differential Equations, 2nd Ed.,  
with contribution of G.N. Milstein and M.B. Nevelson, 
Springer, Berlin, 2012. 

\bibitem{Kovchegov}
Y. Kovchegov, N. Michalowski, 
A Class of Markov Chains with no Spectral Gap, 
Proc. AMS, 141(12) (2013), 4317--4326.

\bibitem{Krylov_selection}
N.V. Krylov, On the selection of a Markov process from a system of processes and the construction of quasi-diffusion processes, 
Math. USSR-Izv., 7(3) (1973), 691--709.


\bibitem{kulik-leonenko}
A.M. Kulik,  N.N. Leonenko, Ergodicity and mixing bounds for the Fisher--Snedecor diffusion. Bernoulli 19(5B), (2013),  2294--2329. 


\bibitem{Eva}
E. L\"{o}cherbach, D. Loukianova, O.  Loukianov, 
Polynomial bounds in the Ergodic theorem for one-dimensional diffusions and integrability of hitting times, 
Ann. Inst. Henri Poincare -- prob. 47 (2011), 425--449. 


\bibitem{MaRo}
Z.-M. Ma, M. R\"ockner, 
Introduction to the Theory of (Non-Symmetric) Dirichlet Forms, 1st ed., Springer, 1992. 

\bibitem{Mao}
Y. Mao, Convergence rates in strong ergodicity for Markov processes, Stoch. Proc. Appl. 116(2006), 1946-1976.


\bibitem{McKean}
H.P. McKean, Stochastic integrals, 
AMS Chelsea Publishing, 1969 (reprinted 2005).


\bibitem{MenshPopov}
M.V. Menshikov, S.Y. Popov, 
Exact power estimates for countable Markov chains, 
Markov Process. Related Fields, 1 (1) (1995) 57--78.

\bibitem{MenshPopov2}
M. Menshikov, S. Popov, A. Wade, 
Non-homogeneous Random Walks: Lyapunov Function Methods for Near-Critical Stochastic Systems, Cambridge, CUP, 2016. 


\bibitem{ayv_grad_drift2}
A. Uglov, A. Veretennikov, 
Yet again on polynomial convergence for SDEs with a gradient-type drift, 
arXiv:1706.09374. 


\bibitem{Ve81}
A.Yu. Veretennikov,  On strong and weak solutions of one - dimensional stochastic equations with boundary conditions. Theory Probab. Appl. 26(4) (1981), 670-686.

\bibitem{ayv_grad_drift}
A.Yu. Veretennikov, On polynomial mixing for SDEs with a gradient-type drift, Theory Probab. Appl., 45(1) (2001) 160--164.



\bibitem{2006}
A.Yu. Veretennikov, 
On lower bounds for mixing coefficients of Markov diffusions, In: From Stochastic Calculus to Mathematical Finance; The Shiryaev Festschrift. Kabanov, Yu.; Liptser, R.; Stoyanov, J. (Eds.) Springer, Berlin et al. (2006), 623--633. 

\bibitem{Ve1} 
A.Yu. Veretennikov,  On polynomial mixing bounds for stochastic
differential equations. Stochastic Processes and their Applications, 70 (1997)
115--127.

\bibitem{Ve2} 
A.Yu. Veretennikov,  On polynomial mixing and convergence rate for
stochastic difference and differential equations, Teoriya Veroyatnostej i ee
Primenen., 44(2) (1999) 312--327.

\bibitem{Ve_LN}
A.Yu. Veretennikov. Ergodic Markov processes and Poisson equations (lecture notes). In: Modern problems of stochastic analysis and statistics - Selected contributions in honor of Valentin Konakov (editor: V. Panov), Springer, 2017, 457 -- 511. 

\bibitem{Zv}
G.A. Zverkina, 
On strong bounds of rate of convergence for regenerative processes, 
Cham, Springer, 2016, 381-393.



\end{thebibliography}

\end{document}